\documentclass[11pt]{article}
\usepackage{amsthm,epsfig,a4}
\usepackage[ansinew]{inputenc}
\usepackage{amsmath,amssymb}
\usepackage{amsfonts}
\usepackage{eufrak}
\def\N{\mathbb{N}}
\def\Z{\mathbb{Z}}
\def\R{\mathbb{R}}
\def\C{\mathbb{C}}
\def\e{\mathrm{e}}
\def\i{\mathrm{i}}
\def\dr{\mathrm{d}}
\def\q{\nolinebreak \hfill $\Box$} 
\def\A{\\[5mm] \indent} 

\def\:{\colon\thinspace}
\theoremstyle{plain}

\newtheorem{SLemma}{Lemma}
\newtheorem{SProp}[SLemma]{Proposition}
\newtheorem{SThm}[SLemma]{Theorem}

\theoremstyle{definition}

\newtheorem*{Rem}{Remark}

\allowdisplaybreaks
\begin{document}
%
%
\title{Odd-dimensional solvmanifolds are contact}
\author{\textsc{Christoph Bock}
}
\date{}
\maketitle


\begin{abstract}
Bourgeois proved in \cite{Bou} that odd-dimensional tori admit a contact structure.
We shall prove a more general result:
Any odd-dimensional parallelisable closed manifold admits a contact structure.
This implies that a solvmanifold $\Gamma \backslash G$ is contact, where $\Gamma$ is a lattice in a connected and simply-connected solvable Lie group $G$ of odd dimension.
\end{abstract}

Throughout this short note a manifold is assumed to be smooth, connected, oriented and to have no boundary.

A \emph{contact manifold} is a pair $(M, \xi)$, where $M$ is a $(2n+1)$-manifold, $n \in \N$, and $\xi$ can be written as $\xi = \ker \alpha$, where $\alpha \in \Omega^1(M)$ is a differential $1$-form with $\alpha_p \wedge {(\dr \alpha)^n}_p \ne 0$ for all $p \in M$.
This implies that the structure group of the tangent bundle of $M$ reduces to ${\mathrm U}(n) \times \{1\}$, see \cite[p.\ 68]{GeB2}.
Therefore, it makes sense to call a manifold \emph{almost contact}, if the structure group of the tangent bundle reduces to ${\mathrm U}(n) \times \{1\}$.
There is the following important consequence of the deep result \cite[Corollary 1.3]{BER} of Borman, Eliashberg and Murphy:

\begin{SProp} \label{almost contact contact}
Let $M$ be an almost contact closed manifold.
Then $M$ admits the structure of a contact manifold. \q
\end{SProp}

Therefore, in order to obtain the next theorem, it is enough to prove the next proposition.

\begin{SThm} \label{Hauptsatz}
Any odd-dimensional parallelisable closed manifold admits a contact structure.
\end{SThm}

\begin{SProp} \label{Hauptsatz1}
Any odd-dimensional parallelisable closed manifold is almost contact.
\end{SProp}

\textit{Proof.} Let $n \in \N$, $M$ be a $(2n+1)$-dimensional closed manifold equipped with a parallelisation $(X_1, \ldots, X_{2n+1})$.
Then,
$$ J X_1 := X_2, \, J X_2 := -X_1, \, \ldots, \, J X_{2n-1} := X_{2n}, \, J X_{2n} := -X_{2n-1} $$
gives rise to an almost contact structure $( \langle X_1, \ldots, X_{2n} \rangle_{\R}, J, \langle X_{2n+1} \rangle_{\R})$ on $M$, see again \cite[p.\ 68]{GeB2}, and the structure group of the tangent bundle of $M$ reduces to ${\mathrm U}(n) \times \{1\}$. \q
\A
In \cite{ipse3}, we asked whether every five-dimensional solvmanifold admits a contact structure.\footnote{There, by another approach, we saw that most five-dimensinal solvmanifolds are contact.}
Using Theorem \ref{Hauptsatz}, we can give the answer for any odd-dimen\-sional solvmanifold.

A \emph{solvmanifold} is a homogeneous space $\Gamma \backslash G$, where $G$ is a connected and simply-connected solvable Lie group and $\Gamma$ a lattice in $G$.
Recall that a \emph{lattice in a Lie group} is a discrete co-compact subgroup.
Therefore, tori are solvmanifolds.

\begin{Rem}
It is important to note that there is a more general notion of solvmanifold, namely a compact quotient of a connected and simply-connected solvable Lie group by a (possibly non-discrete) closed Lie subgroup (see \cite{Aus}), but we are only considering solvmanifolds as in the last definition.
Sometimes, such are called \emph{special solvmanifolds} in the literature.

By \cite[Theorem 2.3.11]{TO}, a solvmanifold in our sense is necessary parallelisable.
E.g.\ the Klein bottle (which can be written as the compact homogeneous space
$$ \{ (\frac{k}{2}, l + \i y) \, | \, k,l \in \Z \, \wedge \, y \in \R \} \backslash (\R \ltimes_{\mu} \C), \mbox{ where } \forall_{t \in \R} \forall_{z \in \C} \, \mu(t)(z) := \e^{2 \pi \i t} z, $$
of the three-dimensional connected and simply-connected solvable Lie group $\R \ltimes_{\mu} \C$) is not a solvmanifold covered by our definition.
(The manifold $M$ in \cite[Section 3]{AS} is not a solvmanifod covered by our definition, too.)
\end{Rem}

\cite[Theorem 2.3.11]{TO} and Theorem \ref{Hauptsatz} imply:

\begin{SThm} \label{Hauptsatz2}
Odd-dimensional solvmanifolds admit a contact structure. \q
\end{SThm}

%
%
%
%

\noindent
\textsc{Christoph Bock\\ Department Mathematik\\ Universität Erlangen-N\"urn\-berg\\ Cauerstraße 11\\ 91058 Erlangen\\ Germany}

\noindent
\textit{e-mail:} \verb"bock@mi.uni-erlangen.de"
\end{document}